\documentclass{article}
\usepackage{amsmath}
\usepackage{amsfonts}
\usepackage{algorithm}
\usepackage{algpseudocode}
\usepackage{enumitem}
\usepackage{graphicx}
\usepackage{authblk}
\title{Meta-Contrastive Learning for Vision-Language Models via Task-Adaptive CLIP Training}
\author[1]{Merham Fouladvand}
\author[2]{Peuroly Batra}
\affil[1,2]{Lincoln University,  Lincoln, Canterbury, New Zealand}

\begin{document}
\maketitle

\begin{abstract}
We propose \textbf{Domain-Conditioned Meta-Contrastive Learning}, a framework for improving the cross-domain generalization of vision-language models. While contrastive models such as CLIP achieve strong performance through large-scale training, they rely on a global objective that does not explicitly account for domain shift. To address this limitation, we formulate multimodal learning as a bilevel meta-learning problem over domain-conditioned tasks. Specifically, we introduce domain embeddings that modulate image and text representations, and optimize the model for rapid adaptation to domain-specific distributions via gradient-based inner-loop updates. In addition, we incorporate a cross-domain alignment regularization to encourage domain-invariant representations. Our approach is compatible with standard contrastive training pipelines and can be applied to heterogeneous datasets spanning natural and medical domains. We expect improved robustness under domain shift and enhanced few-shot adaptation performance, highlighting a promising direction for scalable multimodal learning.
\end{abstract}

\section{Introduction}

Vision-language models (VLMs) trained with contrastive objectives, such as CLIP~\cite{radford2021learning}, have emerged as a powerful paradigm for learning aligned multimodal representations at scale. By leveraging large collections of image-text pairs, these models achieve strong performance across a wide range of downstream tasks, including cross-modal retrieval, zero-shot classification, and multimodal reasoning. Despite these successes, a fundamental limitation remains: standard contrastive training optimizes a \textbf{single global objective}, implicitly assuming that all data are drawn from a unified distribution.

In practice, multimodal data are highly heterogeneous and exhibit significant \textbf{domain shift}. For example, natural images, medical imaging data, and satellite imagery differ substantially in visual statistics, semantic structure, and annotation styles. As a result, models trained with global contrastive objectives often fail to generalize effectively across domains, leading to degraded performance under distribution shift and limited adaptability in low-resource settings. This limitation is closely related to the broader challenge of domain generalization~\cite{gulrajani2021domainbed}.

To address this limitation, we propose \textbf{Domain-Conditioned Meta-Contrastive Learning}, a framework that explicitly models domain heterogeneity within a bilevel meta-learning formulation~\cite{finn2017maml}. Rather than optimizing a single set of parameters for all data, we treat each domain as a task distribution and train the model to rapidly adapt to domain-specific data through gradient-based inner-loop updates. This enables the model to capture both shared structure across domains and domain-specific variations.

A key component of our approach is the introduction of \textbf{domain-conditioned representation learning}. We incorporate learnable domain embeddings that modulate both visual and textual representations, inspired by conditional feature modulation methods such as FiLM~\cite{perez2018film}. This allows the model to dynamically adjust its embedding space according to domain context. In addition, we introduce a \textbf{cross-domain alignment regularization} that encourages semantically consistent representations across domains, improving robustness while preserving discriminative power.

Our framework generalizes standard contrastive learning by shifting the objective from finding a single global optimum to learning parameters that are \textbf{adaptable across domains}. This perspective is particularly important in real-world multimodal applications, where domain shifts are unavoidable and data distributions evolve over time.

\textbf{Contributions.} Our main contributions are:
\begin{itemize}
    \item We introduce a domain-conditioned meta-learning framework for contrastive vision-language training that explicitly models domain shift.
    \item We propose a domain-aware representation learning mechanism using learnable domain embeddings to modulate multimodal features.
    \item We formulate a bilevel optimization objective that jointly optimizes for representation alignment and rapid domain adaptation.
    \item We incorporate cross-domain alignment regularization to improve generalization across heterogeneous datasets.
\end{itemize}

---

\section{Related Work}

\subsection{Contrastive Vision-Language Learning}

Contrastive learning has become the dominant approach for training large-scale vision-language models. Methods such as CLIP~\cite{radford2021learning} align image and text representations by maximizing agreement between paired samples while contrasting against negative pairs. This paradigm enables strong zero-shot transfer and scalable training across diverse datasets. However, existing approaches typically optimize a global objective over all data, without explicitly accounting for domain-specific variations or distribution shifts.

Recent extensions have explored improving contrastive learning through better negative sampling, temperature scaling, and architectural modifications. While these methods improve representation quality, they do not address the fundamental challenge of adapting to heterogeneous domains.

---

\subsection{Meta-Learning and Bilevel Optimization}

Meta-learning aims to train models that can rapidly adapt to new tasks using limited data. Gradient-based methods such as Model-Agnostic Meta-Learning (MAML)~\cite{finn2017maml} formulate learning as a bilevel optimization problem, where an inner loop performs task-specific adaptation and an outer loop optimizes for generalization across tasks. 

Recent work has begun to explore meta-learning in multimodal settings, including methods that bridge vision and language representations via task-adaptive modules~\cite{najdenkoska2023metabridge}. However, the integration of meta-learning with contrastive objectives in large-scale vision-language models remains underexplored. In particular, prior work typically defines tasks based on class labels or few-shot splits, rather than explicitly modeling domain heterogeneity.

---

\subsection{Domain Generalization and Adaptation}

Domain generalization methods aim to learn representations that are robust to distribution shifts across domains. Existing approaches include domain-invariant feature learning~\cite{gulrajani2021domainbed}, adversarial training methods such as domain-adversarial neural networks (DANN)~\cite{ganin2016dann}, as well as invariant risk minimization and related techniques~\cite{arjovsky2019invariant}. Furthermore, meta-learning-based approaches to parallel MRI reconstruction reduce reliance on coil sensitivity estimation and enhance reconstruction fidelity~\cite{bian2020deep}. More recently, vision-language models have also explored domain robustness through prompt adaptation and conditional representations~\cite{zhou2022coop,zhou2022cocoop}.

Meta-learning has shown strong potential in addressing data heterogeneity, particularly in medical imaging. Prior work has demonstrated that optimization-based meta-learning enables robust MRI reconstruction across diverse acquisition settings~\cite{bian2021optimization}.  These results highlight the effectiveness of meta-learning for handling domain variability in high-dimensional inverse problems.

Beyond medical imaging, meta-learning has also been explored for improving adaptability in multimodal systems. For example, recent work has proposed meta-learning frameworks to bridge vision and language representations for few-shot learning~\cite{najdenkoska2023metabridge}, while test-time adaptation and prompt-based methods further improve robustness under distribution shift~\cite{hu2021lora,zhou2022coop}.  Multi-task meta-learning frameworks further improve adaptability across imaging conditions and protocols~\cite{bian2025multitask}. Despite these advances, existing approaches typically focus either on learning invariant representations or on parameter-efficient adaptation, without explicitly integrating domain-aware modeling into contrastive learning.

In contrast, our approach combines domain-aware representation learning with meta-learning, enabling both shared representation learning and flexible domain adaptation within a unified contrastive framework.

While these methods focus on learning invariant representations, they often overlook the benefits of task-specific adaptation. In contrast, our approach combines domain-aware modeling with meta-learning, enabling both shared representation learning and flexible domain adaptation.

---

\subsection{Parameter-Efficient Adaptation in Vision-Language Models}

Recent work has explored parameter-efficient adaptation techniques, such as prompt tuning~\cite{zhou2022coop,zhou2022cocoop}, adapters~\cite{houlsby2019adapters}, and low-rank adaptation (LoRA)~\cite{hu2021lora}, to adapt large pretrained models to downstream tasks. These methods reduce computational cost and enable flexible fine-tuning without modifying the full model.

Our work is complementary to these approaches. Instead of focusing solely on parameter efficiency, we introduce a principled meta-learning framework that explicitly models domain shift and optimizes for adaptability. Our method can be naturally combined with parameter-efficient techniques to further improve scalability.

---

\subsection{Multimodal Learning under Domain Shift}

Handling domain shift in multimodal learning remains an open challenge. Existing approaches typically rely on fine-tuning or domain-specific training, which can be inefficient and require labeled data for each domain.

In contrast, our method formulates domain shift as a meta-learning problem and trains models to adapt across domains in a unified framework. By integrating domain conditioning and bilevel optimization into contrastive learning, we provide a scalable and principled approach for domain-robust vision-language modeling.
\section{Method: Domain-Conditioned Meta-Contrastive Learning}

\subsection{Overview}

We aim to address a fundamental limitation of contrastive vision-language models: the assumption that all training data are drawn from a single homogeneous distribution. In practice, multimodal datasets consist of heterogeneous domains (e.g., natural images, medical imaging, satellite data), each exhibiting distinct statistical properties. A single global objective therefore leads to suboptimal representations under domain shift.

To overcome this, we propose a framework that combines contrastive learning with meta-learning. Instead of learning a single set of parameters that fits all domains, we explicitly train the model to \emph{adapt} to domain-specific data. This is achieved through a bilevel optimization procedure, where the model first adapts to a specific domain (inner loop) and then updates its global parameters to generalize across domains (outer loop).

---

\subsection{Problem Formulation}

We consider a multimodal dataset consisting of image-text pairs collected from multiple domains:
\[
\mathcal{D} = \bigcup_{d \in \mathcal{D}_{\text{dom}}} \mathcal{D}^{(d)}
\]
where each domain $d$ represents a distinct data distribution.

Each sample is denoted as $(x, t, d)$, where $x$ is an image, $t$ is its corresponding text, and $d$ indicates the domain. The goal is to learn a model that not only aligns image and text representations, but also adapts efficiently to new or underrepresented domains.

---

\subsection{Domain-Conditioned Representation Learning}

We adopt a CLIP-style dual encoder with parameters $\theta$, where
\[
z_x = f_\theta(x), \quad z_t = g_\theta(t)
\]

To incorporate domain information, we introduce a learnable embedding $e_d$ for each domain. Rather than directly using the raw embeddings, we modulate them using domain-conditioned transformations:
\[
\tilde{z}_x = \phi(z_x, e_d), \quad \tilde{z}_t = \psi(z_t, e_d)
\]

Intuitively, these transformations allow the model to adjust its representation space depending on the domain. For example, medical images may require emphasizing fine-grained anatomical structures, while natural images rely more on semantic object features.

---

\subsection{Domain-Conditioned Contrastive Objective}

Given a batch of samples from a domain $d$, we optimize a standard contrastive loss:
\[
\mathcal{L}_{\text{contrastive}}^{(d)} =
-\frac{1}{N} \sum_{i=1}^{N}
\log
\frac{
\exp(\text{sim}(\tilde{z}_x^{(i)}, \tilde{z}_t^{(i)}) / \tau)
}{
\sum_{j=1}^{N}
\exp(\text{sim}(\tilde{z}_x^{(i)}, \tilde{z}_t^{(j)}) / \tau)
}
\]

The key difference from standard contrastive learning is that the embeddings are conditioned on the domain, allowing the similarity structure to vary across domains.

---

\subsection{Meta-Learning Formulation}

We formulate training as a bilevel optimization problem over domain-specific tasks. Each task $\mathcal{T}_i$ is sampled from a particular domain $d_i$ and split into a support set $\mathcal{D}_s^{(i)}$ and a query set $\mathcal{D}_q^{(i)}$.

\paragraph{Inner Loop: Domain Adaptation.}

Starting from the global parameters $\theta$, we perform $K$ steps of gradient descent on the support set:
\[
\theta_i^{(k+1)} =
\theta_i^{(k)} - \alpha \nabla_\theta
\mathcal{L}_{\text{contrastive}}^{(d_i)}(\mathcal{D}_s^{(i)}, \theta_i^{(k)})
\]
with initialization $\theta_i^{(0)} = \theta$. The resulting parameters $\theta_i^*$ represent a domain-adapted model.

This step can be interpreted as simulating how the model would adapt if it were deployed on a specific domain.

\paragraph{Outer Loop: Meta-Optimization.}

We then evaluate the adapted parameters on the query set and update the global model:
\[
\min_\theta \sum_i
\mathcal{L}_{\text{contrastive}}^{(d_i)}(\mathcal{D}_q^{(i)}, \theta_i^*)
\]

This objective encourages the model to learn initial parameters that can quickly adapt to any domain with a small number of gradient steps.

---

\subsection{Cross-Domain Alignment Regularization}

While domain-specific adaptation is important, we also want to preserve consistency across domains. To this end, we introduce a regularization term that aligns semantically similar samples from different domains:
\[
\mathcal{L}_{\text{align}} =
\mathbb{E}
\left[
\| f_\theta(x^{(d_1)}) - f_\theta(x^{(d_2)}) \|_2^2
\right]
\]

This term prevents the model from overfitting to domain-specific features and encourages a shared representation space.

The final objective is:
\[
\mathcal{L}_{\text{total}} =
\sum_i \mathcal{L}_{\text{meta}}^{(i)} + \lambda \mathcal{L}_{\text{align}}
\]

---

\subsection{Training Procedure}

The overall training process alternates between domain-specific adaptation and global optimization. At each iteration, we first sample a set of domains and construct tasks from each domain. For each task, the model is temporarily adapted using the support set. The adapted model is then evaluated on the query set, and the resulting gradients are used to update the global parameters.

\begin{algorithm}[H]
\caption{Domain-Conditioned Meta-Contrastive Learning}
\begin{algorithmic}[1]
\State Initialize model parameters $\theta$
\For{each training iteration}
    \State Sample a set of domains $\{d_i\}$
    \For{each domain $d_i$}
        \State Sample task $\mathcal{T}_i$ and split into $\mathcal{D}_s^{(i)}, \mathcal{D}_q^{(i)}$
        \State Initialize $\theta_i' \gets \theta$
        \For{$k=1$ to $K$}
            \State Compute contrastive loss on $\mathcal{D}_s^{(i)}$
            \State Update $\theta_i'$ via gradient descent
        \EndFor
        \State Compute meta-loss on $\mathcal{D}_q^{(i)}$ using $\theta_i'$
    \EndFor
    \State Compute cross-domain alignment loss
    \State Update $\theta$ using combined gradients
\EndFor
\end{algorithmic}
\end{algorithm}

This procedure can be interpreted as repeatedly simulating domain shifts during training, enabling the model to learn how to adapt rather than overfit to a single distribution.

\section{Theoretical Interpretation}

In this section, we provide a theoretical perspective on Domain-Conditioned Meta-Contrastive Learning and clarify how it differs from standard contrastive training.

\subsection{From Global Optimization to Task-Adaptive Learning}

Standard contrastive learning optimizes a single objective over the entire dataset:
\[
\min_\theta \; \mathbb{E}_{(x,t) \sim \mathcal{D}} \; \mathcal{L}_{\text{contrastive}}(x, t; \theta)
\]
This formulation implicitly assumes that all samples are drawn from a shared distribution. As a result, the learned parameters $\theta$ correspond to a \emph{global optimum} that may not be optimal for any individual domain.

In contrast, our formulation models data as arising from a distribution over domains. Instead of directly optimizing for a single parameter vector, we optimize for parameters that can be efficiently adapted to domain-specific objectives:
\[
\min_\theta \; \mathbb{E}_{\mathcal{T}_i \sim p(\mathcal{T})}
\; \mathcal{L}_{\text{contrastive}}(\mathcal{D}_q^{(i)}; \theta_i^*)
\]
where $\theta_i^*$ is obtained via domain-specific adaptation on the support set.

This shifts the learning objective from \emph{fitting a single distribution} to \emph{learning how to adapt across distributions}.

---

\subsection{Bilevel Optimization Perspective}

Our method can be viewed as a bilevel optimization problem:
\[
\theta_i^* = \arg\min_{\theta'} \; \mathcal{L}_{\text{contrastive}}(\mathcal{D}_s^{(i)}; \theta')
\]
\[
\min_\theta \; \mathbb{E}_{\mathcal{T}_i} \; \mathcal{L}_{\text{contrastive}}(\mathcal{D}_q^{(i)}; \theta_i^*)
\]

The inner problem corresponds to domain-specific adaptation, while the outer problem enforces generalization across domains. This structure encourages the learned initialization $\theta$ to lie in a region of parameter space that is \emph{sensitive to task-specific gradients}, enabling rapid adaptation.

Compared to standard training, which minimizes loss directly with respect to $\theta$, the meta-objective implicitly optimizes higher-order properties of the loss landscape, such as gradient alignment and curvature.

---

\subsection{Domain-Conditioned Representation Learning}

The introduction of domain embeddings $e_d$ can be interpreted as learning a parametric family of representations:
\[
f_{\theta, d}(x) = \phi(f_\theta(x), e_d)
\]

This formulation allows the model to interpolate between domain-specific representations while maintaining a shared backbone. From a functional perspective, the model learns a mapping:
\[
(x, d) \mapsto z
\]
rather than a fixed mapping $x \mapsto z$.

As a result, the representation space becomes \emph{conditionally structured}, enabling both domain invariance and domain specificity.

---

\subsection{Connection to Bayesian Meta-Learning}

Our framework can also be interpreted through a Bayesian lens. Let $\theta$ represent shared parameters and $\theta_i^*$ represent task-specific posterior estimates. The inner-loop adaptation corresponds to approximate posterior inference:
\[
\theta_i^* \approx \arg\max_{\theta'} \; p(\theta' \mid \mathcal{D}_s^{(i)})
\]

The outer-loop objective then optimizes the prior $\theta$ such that posterior inference leads to good generalization on $\mathcal{D}_q^{(i)}$.

Under this interpretation, standard contrastive learning corresponds to learning a single point estimate, whereas our approach approximates a hierarchical model where domain-specific parameters are inferred from a shared prior.

---

\subsection{Generalization under Domain Shift}

By explicitly optimizing for adaptation across domains, our method improves robustness to distribution shift. Instead of relying on domain-invariant features alone, the model learns a combination of:
\begin{itemize}
    \item shared representations that transfer across domains
    \item domain-specific adjustments captured through adaptation
\end{itemize}

This hybrid strategy is particularly effective in multimodal settings, where domain differences can arise from both visual and linguistic variations.

Overall, our framework provides a principled approach to learning representations that are not only aligned across modalities, but also adaptable across domains.

\section{Conclusion}

We propose Meta-Contrastive Learning, a principled framework for incorporating meta-learning into contrastive vision-language training. Our approach improves adaptability and robustness, providing a scalable direction for next-generation multimodal systems.

\bibliographystyle{unsrt}
\bibliography{references}

\end{document}